\documentclass{article}

\usepackage{hyperref}

\begin{document}

\noindent
{\sc published elec.$\,$at}

\vspace*{-0.5 mm}

\noindent
\href{http://www.maplesoft.com/}{\texttt{www.maplesoft.com}}

\medskip

\medskip
\begin{center}
{\Large \boldmath \bf
A PROCEDURE FOR FINDING                                         \\[1.0 ex]
THE $k^{\mbox{\footnotesize \bf TH}}$ POWER OF A MATRIX}
\end{center}

\bigskip
\begin{center}
{\large \it Branko Male\v sevi\' c}                             \\
Faculty of Electrical Engineering,                              \\
University of Belgrade, Serbia                                  \\
\href{mailto:malesh@EUnet.yu}{\textrm{malesh@EUnet.yu}}
\end{center}

\smallskip
\begin{center}
{\large \it Ivana Jovovi\' c}                                   \\
PhD student, Faculty of Mathematics,                            \\
University of Belgrade, Serbia                                  \\
\href{mailto:ivana121@EUnet.yu}{\textrm{ivana121@EUnet.yu}}
\end{center}

\medskip

\section{Introduction }

\smallskip

\medskip
\noindent
This worksheet demonstrates the use of Maple in Linear Algebra.

\medskip
\noindent
We give a new procedure ({\tt PowerMatrix}) in Maple for finding the $k^{\mbox{\footnotesize \rm th}}$
power of $n$-by-$n$ square matrix $\mbox{\tt A}$, in a symbolic form, for any positive integer $k$, $k \!\geq\! n$.
The algorithm is based on an application of Cayley-Hamilton theorem. We used the fact that the entries of the matrix
$\mbox{\tt A}^k$ satisfy the same recurrence relation which is determined by the characteristic polynomial of
the matrix $\mbox{\tt A}$ (see~[{\bf 1}]). The order of these recurrences is $n-d$, where $d$ is the lowest degree of
the characteristic polynomial of the matrix $\mbox{\tt A}$.

\medskip
\noindent
For non-singular matrices the procedure can be extended for $k$ not only a positive integer.

\smallskip

\section{Initialization}

\quad$\;\!\!\!\!\!\!\!\!>\;$ {\tt restart:}

\noindent
\quad$\;${\tt with(LinearAlgebra):}

\medskip

\subsection{Procedure Definition}

\medskip

\subsubsection{PowerMatrix}

\noindent
Input data are a square matrix {\tt A} and a parameter $k$. Elements of  the matrix {\tt A}
can be numbers and/or parameters. The parameter $k$ can take numeric value or be a symbol.
The output data is the $k^{\mbox{\footnotesize \rm th}}$ power of the matrix. The procedure
{\tt PowerMatrix} is as powerful as the procedure {\tt rsolve}.

\smallskip
$\!\!\!\!\!\!\!\!>\;$
{\tt PowerMatrix$\,:=\;$proc(A::Matrix,k)}

\smallskip
{\tt local i,j,m,r,q,n,d,f,P,F,C;}

\smallskip
{\tt P$\,:=\;$x->CharacteristicPolynomial(A,x);}

\smallskip
{\tt n$\,:=\;$degree(P(x),x);}

\smallskip
{\tt d$\,:=\;$ldegree(P(x),x);}

\break

\smallskip
{\tt F$\,:=\;$(i,j)->rsolve(sum(coeff(P(x),x,m)*f(m+q),m=0..n)=0,seq(f(r)=(A\verb|^|r)[i,j],}

\smallskip
{\tt r=d+1..n),f);}

\smallskip
{\tt C$\,:=\;$q->Matrix(n,n,F);}

\smallskip
{\tt if (type(k,integer)) then return(simplify(A\verb|^|k)) elif (Determinant(A)=0 and }

\smallskip
{\tt not type(k,numeric)) then printf("The \%a${}^{\footnotesize \tt th}$ power of the matrix for
\%a>=\%d:"$\!$,}

\smallskip
{\tt k,k,n) elif (Determinant(A)=0 and type(k,numeric)) then return(simplify(A\verb|^|k)) fi;}

\smallskip
{\tt return(simplify(subs(q=k,C(q))));}

\smallskip
{\tt end:}

\section{Examples}

\medskip

\subsection{Example 1.}

\smallskip
\quad$\,$
$\!\!\!\!\!\!\!\!>\;$
{\tt A$\,:=\;$Matrix([[4,-2,2],[-5,7,-5],[-6,6,-4]]);}

$$
\mbox{\tt A$\,:=$}\left[
{\begin{array}{rrr}
4 & -2 & 2 \\
-5 & 7 & -5 \\
-6 & 6 & -4
\end{array}}
 \right]
$$

\smallskip
$\!\!\!\!\!\!\!\!>\;$
{\tt PowerMatrix(A,k);}

$$
\left[
{\begin{array}{ccc}
 - 2^{k} + 2 \cdot 3^{k} & 2^{(1 + k)} - 2 \cdot 3^{k} &  - 2^{(1 + k)} + 2 \cdot 3^{k}                 \\
 - 5 \cdot 3^{k} + 5 \cdot 2^{k} & 5 \cdot 3^{k} - 4 \cdot 2^{k} &  - 5 \cdot 3^{k} + 5 \cdot 2^{k}     \\
6 \cdot 2^{k} - 6 \cdot 3^{k} &  - 6 \cdot 2^{k} + 6 \cdot 3^{k} &  - 6 \cdot 3^{k} + 7 \cdot 2^{k}
\end{array}}
\right]
$$

\smallskip
$\!\!\!\!\!\!\!\!>\;$
{\tt Determinant(A);}{

$$
12
$$

$\!\!\!\!\!\!\!\!>\;$
{\tt B$\,:=\;$A\verb|^|(-1);}

$$
\mbox{\tt B$\,:=$}\left[
{\begin{array}{rrr}
 \mbox{\small $\displaystyle\frac{1}{6}$} &  \mbox{\small $\displaystyle \frac{1}{3}$}  &
-\mbox{\small $\displaystyle\frac{1}{3}$}                                        \\ [2ex]
 \mbox{\small $\displaystyle\frac{5}{6}$} & -\mbox{\small $\displaystyle \frac{1}{3}$}  &
 \mbox{\small $\displaystyle \frac{5}{6}$}                                       \\ [2ex]
1 & -1 & \mbox{\small $\displaystyle \frac{3}{2}$}
\end{array}}
\right]
$$

\smallskip
$\!\!\!\!\!\!\!\!>\;$
{\tt PowerMatrix(B,k);}

$$
\left[
{\begin{array}{ccc}
 - 2^{( - k)} + 2 \cdot 3^{( - k)} & 2^{(1 - k)} - 2 \cdot 3^{( - k)} &  - 2^{(1 - k)} + 2 \cdot 3^{( - k)} \\
 - 5 \cdot 3^{( - k)} + 5 \cdot 2^{( - k)} & 5 \cdot 3^{( - k)} - 4 \cdot 2^{( - k)}  & - 5 \cdot 3^{( - k)} +
   5 \cdot 2^{( - k)} \\
 - 6 \cdot 3^{( - k)} + 6 \cdot 2^{( - k)} &  - 6 \cdot 2^{( - k)} + 6 \cdot 3^{( - k)} &
 - 6 \cdot 3^{( - k)} + 7 \cdot 2^{( - k)}
\end{array}}
\right]
$$

\medskip

\subsection{Example 2.}

\smallskip
\quad$\,$
$\!\!\!\!\!\!\!\!>\;$
{\tt A$\,:=\;$Matrix([[1-p,p],[p,1-p]]);}

$$
\mbox{\tt A$\,:=$}\left[
{\begin{array}{cc}
1 - p & p \\
p & 1 - p
\end{array}}
\right]
$$

\break

\smallskip
$\!\!\!\!\!\!\!\!>\;$
{\tt PowerMatrix(A,k);}

$$
\left[
{\begin{array}{cc}
{\displaystyle \frac {(1 - 2\,p)^{k}}{2}}  + {\displaystyle
\frac {1}{2}}  &  - {\displaystyle \frac {(1 - 2\,p)^{k}}{2}}  +
{\displaystyle \frac {1}{2}}  \\ [2ex]
 - {\displaystyle \frac {(1 - 2\,p)^{k}}{2}}  + {\displaystyle
\frac {1}{2}}  & {\displaystyle \frac {(1 - 2\,p)^{k}}{2}}  +
{\displaystyle \frac {1}{2}}
\end{array}}
\right]
$$

\medskip
{\tt The example is from [{\bf 4}], page 272, exercise 19.}

\medskip

\subsection{Example 3.}

\quad$\,$
$\!\!\!\!\!\!\!\!>\;$
{\tt A$\,:=\;$Matrix([[a,b,c],[d,e,f],[g,h,i]]);}

$$
\mbox{\tt A$\,:=$}\left[
{\begin{array}{ccc}
a & b & c \\
d & e & f \\
g & h & i
\end{array}}
\right]
$$

$\!\!\!\!\!\!\!\!>\;$
{\tt PowerMatrix(A,k)[1,1];}

$$
\begin{array}{rl}
  & \sum_{\mbox{\scriptsize $\underline{\;}\mbox{\rm R} = \mbox{\footnotesize \tt RootOf}\mbox{\small (}\,(
    g b f
  + h d c
  + i e a
  - g c e
  - h f a
  - i d b) \underline{\;}Z^3
  + ( g c
    + h f
    + d b
    - i e
    - i a
    - e a) \underline{\;}Z^2
 + (i+e+a) \underline{\;}Z -1 \, \mbox{\small )}$}}                                                           \\[2.0 ex]
 & \quad
{\Big [}{\big (}  \underline{\,\,}\mbox{\small \rm R}^2 i e
                - \underline{\,\,}\mbox{\small \rm R}^2 h f
                - \underline{\,\,}\mbox{\small \rm R} e
                - \underline{\,\,}\mbox{\small \rm R} i + 1 {\big )}\,
  \mbox{\Big (}\!\;\mbox{\small $\displaystyle\frac{1}{\underline{\,\,}\mbox{\small \rm R}}$}\;\!\mbox{\Big )}^k
  {\Big /}\;{\big (}( 3 \underline{\,\,}\mbox{\small \rm R}^2 g b f
                    + 3 \underline{\,\,}\mbox{\small \rm R}^2 h d c
                    + 3 \underline{\,\,}\mbox{\small \rm R}^2 i e a                                            \\[2.0 ex]
  & \quad
                    - 3 \underline{\,\,}\mbox{\small \rm R}^2 g c e
                    - 3 \underline{\,\,}\mbox{\small \rm R}^2 h f a
                    - 3 \underline{\,\,}\mbox{\small \rm R}^2 i d b
                    + 2 \underline{\,\,}\mbox{\small \rm R} g c
                    + 2 \underline{\,\,}\mbox{\small \rm R} h f
                    + 2 \underline{\,\,}\mbox{\small \rm R} d b
                    - 2 \underline{\,\,}\mbox{\small \rm R} i e
                    - 2 \underline{\,\,}\mbox{\small \rm R} i a
                    - 2 \underline{\,\,}\mbox{\small \rm R} e a
                    + \, i + e + a   ) \underline{\,\,}\mbox{\small \rm R} {\big )}{\Big ]}
\end{array}
$$

\medskip
$\!\!\!\!\!\!\!\!\#\;$
{\tt Warning!}

{\tt In this example MatrixPower and MatrixFuction procedures cannot be done in real$\!$-time.}

\medskip
$\!\!\!\!\!\!\!\!\#\;$
{\tt MatrixPower(A,k)[1,1];}

\medskip
$\!\!\!\!\!\!\!\!\#\;$
{\tt MatrixFunction(A,v\verb|^|k,v)[1,1];}

\medskip

\subsection{Example 4.}

\smallskip
\quad$\,$
$\!\!\!\!\!\!\!\!>\;$
{\tt A$\,:=\;$Matrix([[0,0,1,0,1],[1,0,0,0,1],[0,0,0,1,1],[0,1,0,0,1],[1,1,1,1,0]]);}

$$
\mbox{\tt A$\,:=$}\left[
{\begin{array}{rrrrr}
0 & 0 & 1 & 0 & 1 \\
1 & 0 & 0 & 0 & 1 \\
0 & 0 & 0 & 1 & 1 \\
0 & 1 & 0 & 0 & 1 \\
1 & 1 & 1 & 1 & 0
\end{array}}
\right]
$$

$\!\!\!\!\!\!\!\!>\;$
{\tt PowerMatrix(A,k)[1,5];}

$$
-
{\displaystyle \frac {\sqrt{17}\,{\bigg (}{\displaystyle \frac {1}{2}}
-
{\displaystyle \frac {\sqrt{17}}{2}}{\bigg )}^{k}}{17}}
+
{\displaystyle \frac {\sqrt{17}\,{\bigg (}{\displaystyle \frac {1}{2}}
+
{\displaystyle \frac {\sqrt{17}}{2}}{\bigg )}^{k}}{17}}
$$

\medskip
{\tt Replace ':' with ';' and see result!}

\medskip
$\!\!\!\!\!\!\!\!>\;$
{\tt MatrixPower(A,k)[1,5]:}

\medskip
$\!\!\!\!\!\!\!\!>\;$
{\tt assume(m::integer):simplify(MatrixPower(A,k)[1,5]):}

\medskip
{\tt The example is from [{\bf 3}], page 101.}

\break

\subsection{Example 5. and Example 6.}

\smallskip
{\tt Pay attention what happens for singular matrices.}

\subsubsection{Example 5.  }

\smallskip
\quad$\,$
$\!\!\!\!\!\!\!\!>\;$
{\tt A$\,:=\;$Matrix([[0,2,1,3],[0,0,-2,4],[0,0,0,5],[0,0,0,0]]);}

$$
\mbox{\tt A$\,:=$}\left[
{\begin{array}{rrrr}
0 & 2 &  1 & 3 \\
0 & 0 & -2 & 4 \\
0 & 0 &  0 & 5 \\
0 & 0 &  0 & 0
\end{array}}
\right]
$$

\smallskip
$\!\!\!\!\!\!\!\!>\;$
{\tt PowerMatrix(A,2);}

$$
\left[
{\begin{array}{rrrr}
0 & 0 & -4 &  13 \\
0 & 0 &  0 & -10 \\
0 & 0 &  0 &   0 \\
0 & 0 &  0 &   0
\end{array}}
\right]
$$

\smallskip
$\!\!\!\!\!\!\!\!>\;$
{\tt PowerMatrix(A,3);}

$$
\left[
{\begin{array}{rrrr}
0 & 0 & 0 & -20 \\
0 & 0 & 0 &   0 \\
0 & 0 & 0 &   0 \\
0 & 0 & 0 &   0
\end{array}}
\right]
$$

\smallskip
$\!\!\!\!\!\!\!\!>\;$
{\tt PowerMatrix(A,k);}

\smallskip
{\tt The ${k}^{\footnotesize \tt th}$ power of the matrix {\tt A} for $k \geq \mbox{\tt 4}$:}

$$
\left[
{\begin{array}{rrrr}
0 & 0 & 0 & 0 \\
0 & 0 & 0 & 0 \\
0 & 0 & 0 & 0 \\
0 & 0 & 0 & 0
\end{array}}
\right]
$$

\medskip
$\!\!\!\!\!\!\!\!>\;$
{\tt MatrixPower(A,k);}

\smallskip
{\tt Error, (in LinearAlgebra:-LA\verb|_|Main:-MatrixPower)}

{\tt power $k$ is not defined for this Matrix}

\bigskip
$\!\!\!\!\!\!\!\!>\;$
{\tt  MatrixFunction(A,v\verb|^|k,v);}

\smallskip
{\tt Error, (in LinearAlgebra:-LA\verb|_|Main:-MatrixFunction)}

{\tt Matrix function $\mbox{\tt v}^k$ is not defined for this Matrix}

\bigskip
{\tt The example is from [\mbox{\bf 2}], page 151, exercise 23.}

\subsubsection{Example 6.}

\smallskip
\quad$\,$
$\!\!\!\!\!\!\!\!>\;$
{\tt A$\,:=\;$Matrix([[1,1,1,0],[1,1,1,-1],[0,0,-1,1],[0,0,1,-1]]);}

$$
\mbox{\tt A$\,:=$}\left[
{\begin{array}{rrrr}
1 & 1 &  1 &  0 \\
1 & 1 &  1 & -1 \\
0 & 0 & -1 &  1 \\
0 & 0 &  1 & -1
\end{array}}
\right]
$$

\break

\smallskip
$\!\!\!\!\!\!\!\!>\;$
{\tt PowerMatrix(A,k);}

\smallskip
$\!\!\!\!\!\!\!\!>\;$
{\tt The ${k}^{\footnotesize \tt th}$ power of the matrix for $k \geq \mbox{\tt 4}$:

$$
\left[
{\begin{array}{cccc}
2^{( - 1 + k)} & 2^{( - 1 + k)} & {\displaystyle \frac {(-1)^{(1 + k)} \cdot 2^{k}}{16}}
+ {\displaystyle \frac {5 \cdot 2^{k}}{16}} &  - {\displaystyle \frac {2^{k}}{16}}
+ {\displaystyle \frac {(-1)^{k} \cdot 2^{k}}{16}}                              \\ [2.5 ex]
2^{( - 1 + k)} & 2^{( - 1 + k)} & {\displaystyle \frac {5 \cdot 2^{k} }{16}}
+ {\displaystyle \frac {5 \cdot (-1)^{(1 + k)} \cdot 2^{k}}{16}}
& {\displaystyle \frac {5 \cdot (-1)^{k} \cdot 2^{k}}{16}}
- {\displaystyle \frac {2^{k}}{16}}                                             \\ [3.0 ex]
0 & 0 & (-1)^{k} \cdot 2^{(-1+k)} & (-1)^{(1+k)} \cdot 2^{(-1+k)}               \\ [3.0 ex]
0 & 0 & (-1)^{(1+k)} \cdot 2^{(-1+k)} & (-1)^{k} \cdot 2^{(-1+k)}
\end{array}}
 \right]
$$

\medskip
$\!\!\!\!\!\!\!\!>$
{\tt MatrixPower(A,k);}

\smallskip
{\tt Error, (in LinearAlgebra:-LA\verb|_|Main:-MatrixPower)}

{\tt power $k$ is not defined for this Matrix}

\bigskip
$\!\!\!\!\!\!\!\!>$
{\tt MatrixFunction(A,v\verb|^|k,v);}

\smallskip
{\tt Error, (in LinearAlgebra:-LA\verb|_|Main:-MatrixFunction)}

{\tt Matrix function $\mbox{\tt v}^k$ is not defined for this Matrix}

\section{References}

\noindent
[{\bf 1}] {\rm Branko Male\v sevi\' c}:$\;${\it Some combinatorial aspects of the composition of a set of functions},~{\tt NSJOM}
{\rm 2006$\;$({\bf 36}),$\;$3-9},$\;$URLs:$\;$
\href{http://www.im.ns.ac.yu/NSJOM/Papers/36_1/NSJOM_36_1_003_009.pdf}{\texttt{http:/$\!$/www.im.ns.ac.yu/NSJOM/Papers/36\mbox{\underline{$\,\,$}}1/NSJOM\mbox{\underline{$\,\,$}}36\mbox{\underline{$\,\,$}}1\mbox{\underline{$\,\,$}}003\mbox{\underline{$\,\,$}}009.pdf}},

\noindent
\href{http://arxiv.org/abs/math.CO/0409287}{\texttt{http:/$\!$/arxiv.org/abs/math.CO/0409287}}.

\medskip

\noindent
[{\bf 2}] {\rm John$\,$B.$\,$Johnston, G.$\,$Baley$\,$Price, Fred$\,$S.$\,$Van Vleck}:$\;${\it Linear Equations
and Matrices},$\;${\rm Addi\-son-Wesley, 1966.}

\medskip

\noindent
[{\bf 3}] {\rm Carl D.$\,$Meyer}:$\;${\it Matrix Analysis and Applied Linear Algebra Book and Solutions Manual}~{\tt SIAM},
{\rm 2001.}

\medskip

\noindent
[{\bf 4}] {\rm Robert Messer}:$\;${\it Linear Algebra Gateway to Mathematics}, {\rm New York,
Harper-Collins~Coll\-ege Publisher, 1993.}

\section{Conclusions}

{\rm This procedure has an educational character. It is an interesting demonstration for finding the
$k^{\footnotesize \rm th}$ power of a matrix in a symbolic form. Sometimes, it gives solutions in the better
form than the existing procedure {\tt MatrixPower} (see example 4.). See also example 5. and example 6.,
where we consider singular matrices. In these cases the procedure {\tt MatrixPower} does not give a solution.
The procedure {\tt PowerMatrix} calculates the $k^{\footnotesize \rm th}$ power of any singular matrices.
In some examples it is possible to get a solution in the better form with using the procedure {\tt allvalues}
(see~example~3.).}

\bigskip
\noindent
{\it Legal Notice$:$
The copyright for this application is owned by the authors. Neither Maplesoft nor the author are responsible
for any errors contained within and are not liable for any damages resulting from the use of this material.
This application is intended for non-commercial, non-profit use only. Contact the author for
permission if you wish to use this application in for-profit activities.}

\break

\end{document}